\documentclass[11pt,twoside,a4paper]{article}
\usepackage[T1]{fontenc}
\usepackage{fancyhdr}
\usepackage{amsmath}
\usepackage{amssymb}
\usepackage{amsfonts}
\usepackage{dsfont}
\usepackage{hyperref}
\usepackage{amsthm}
\usepackage{latexsym}
\usepackage{graphicx}
\usepackage{subfigure}
\usepackage{epsfig}

\usepackage{natbib}
\usepackage{layout}

\usepackage{mathrsfs}

\newtheorem{theorem}{Theorem}[section]
\newtheorem{lemma}[theorem]{Lemma}
\newtheorem{e-proposition}[theorem]{Proposition}
\newtheorem{corollary}[theorem]{Corollary}
\newtheorem{e-definition}[theorem]{Definition}
\newtheorem{remark}{\rm Remark\/}

\newtheorem{theoreme}{Th\'{e}or\`{e}me}[section]

\newtheorem{proposition}[theoreme]{Proposition}

\begin{document}

\title{\large \bf Piecewise linear density estimation for sampled data}
\author{Fran\c{c}ois-Xavier Lejeune\footnote{\tiny L.S.T.A., Universit\'e Paris 6, 175, rue du Chevaleret,
Bo\^ite courrier 158, 75013 Paris, France; E-mail: francois-xavier.lejeune@upmc.fr}}

\maketitle
\date

\noindent{\bf Abstract --} Nonparametric density estimation is considered for a discretely observed stationary continuous-time process. For each of three given time sampling procedures either random or deterministic, we establish that histograms and frequency polygons can reach the same optimal
$L_{2}$-rates as in the independent and identically distributed case. Moreover, thanks to a suitable ``high frequency'' sampling design, these rates are derived together with a minimized time of observation depending on the regularity of sample paths.
\medskip

\medskip


\medskip

\section{Introduction} \label{intro}
Consider a $\mathbb{R}^{d}$-valued process $\{X_{t}, t \in \mathbb{R}\}$, $d \geq 1$,
where all $X_{t}$'s have the same unknown marginal density $f$. The aim of this paper is to study the rates of some nonparametric
piecewise linear estimators of $f$ when the process is discretely sampled in time at $t=t_{1},\dots,t_{n}$. During the past three decades, the problem of density estimation for continuous-time observations has been a subject of continued interest in the statistical literature. Especially, it was shown by Castellana and Leadbetter \cite{CastellanaLeadbetter1986} that if a continuous-time process, observed over the time interval $[0,T]$, has enough irregular sample paths, then nonparametric estimators can achieve a mean-square parametric rate of convergence $1/T$. An account of the research in this field may be found e.g. in two complementary monographs by Bosq and Blanke \cite{BosqBlanke2007} and Kutoyants \cite{Kutoyants2003}, and in Lejeune \cite{Lejeune2006} for the particular case of piecewise linear estimators.

In practice, however, the whole sample path is not always perfectly
observable over a given time period -- either due to technical
reasons or unavailability of data at all time points. Indeed, most
of physical phenomenons usually represented by curves generate rather discrete observed values or
interpolated ones. Hence it seems more natural to plan an estimation approach
based upon $n$ discrete values of the process collected with a suitable time sampling procedure.
In the context of nonparametric density estimation, the three sampling procedures considered in the present work have been investigated by Masry \cite{Masry1983}, Prakasa Rao \cite{PrakasaRao1990}, Wu \cite{Wu1997} (random sampling), Bosq \cite{Bosq1995,Bosq1997} and Blanke and Pumo \cite{BlankePumo2003} (``high frequency'' periodic sampling), among others. As far as we know, most of existing papers -- including those cited above -- only deal with kernel estimation, and none have yet focused a special attention on piecewise linear estimators. In that framework, we are then interested in the rates of histogram and frequency polygon estimators with respect to the mean integrated squared error (MISE) criterion. Here, the chosen histogram-based density estimators have the desirable property of being quickly computed and updated. This is therefore a clear advantage for many applications where typically one has to handle large amounts of data in real time. It is noteworthy that elementary estimators may also be efficient from a theoretical viewpoint. Thus and despite
its high simplicity, the frequency polygon -- defined in dimension one as the linear
interpolant of the mid-points of an equally spaced histogram -- is known to be as good as some
more sophisticated density estimators in terms of MISE (see Scott \cite{Scott1985}).

In this paper, we will show that, under mild conditions,
these estimators built with sampled data have at least the same optimal rates $n^{-2/(d+2)}$ (histogram) and $n^{-4/5}$ (univariate frequency polygon) as in the independent and identically distributed (i.i.d.) case. First, we will exami\-ne the case of two classical random sampling designs, which are a relevant way to treat the occurrence of low frequency and irregularly spaced measurements. Next, we will investigate a deterministic design that applies when the data are observed at high frequency and during a long time, as in a variety of domains, including econometrics, meteorology, oceanology and many others. Particularly, this sampling design is well-adapted to the continuous-time context since the optimal rates can be derived together with a minimized time of observation depending on the regularity of sample paths (see Bosq \cite{Bosq1997}). Thanks to this methodology we will furthermore address the important issue of finding an optimal sampling strategy.

The paper is organized as follows. In Section \ref{Preliminaries}
we will review the time sampling procedures and define our framework. Section \ref{histogram}
contains the main assumptions and our results relative to histograms; the behavior of
frequency polygons is then studied in Section \ref{FreqPoly} and a concluding discussion is given in Section \ref{conclusion}. Finally, the proofs are postponed until Section \ref{Proofs}.

\section{Preliminaries and notations} \label{Preliminaries}
Let $X^{T}=\{X_{t}, 0 \leq t \leq T\}$ be a measurable
$\mathbb{R}^{d}$-valued, $d \geq 1$, continuous-time process on
the probability space $(\Omega, \mathcal{F}, P)$, where the
$X_{t}$'s have a common distribution admitting a density $f$
with respect to the Lebesgue measure over
$\mathbb{R}^{d}$. We suppose that the joint density
$f_{(X_{s},X_{t})}$ of $(X_{s},X_{t})$ does exist for all $s \neq
t$ such that $f_{(X_{s},X_{t})}=f_{(X_{0},X_{|t-s|})}=:f_{|t-s|}$,
which is a quite weak stationary condition (see e.g. Bosq \cite{Bosq1998}). We also denote by
$g_{u}$ the function defined for all $u>0$ as $g_{u}:=f_{u}-f
\otimes f$ where $(f \otimes f)(y,z)=f(y)f(z)$. Some required
asymptotic independence conditions on the process (including
$\alpha$-mixing condition) will be added later with our
assumptions. Our purpose is to estimate the function $f$ from $n$
observations collected up to time $T$ by making use of one of the sampling
procedures described below.

\subsection{Sampling schemes}
Let $\mathscr{T}^{n}=\{t_{k}, 0 \leq k \leq n\}$ be a strictly
increasing sequence of points in time -- or \emph{event arrival times} --
such that $0 = t_{0} < t_{1} < \cdots < t_{n} =: T_{n}$ and $T_{n}
\rightarrow \infty$ as $n \rightarrow \infty$. If $\mathscr{T}^{n}$ is random, it is also
assumed that the processes $X^{T}$ and $\mathscr{T}^{n}$ are
independent and that all $X_{t_{k}}$'s are measurable
with respect to the $\sigma$-algebra generated by $X^{T}$
and $\mathscr{T}^{n}$. The first two random schemes as defined in Masry \cite{Masry1983} are the
following. \vspace{.3cm}

\noindent{\textbf{Renewal sampling}} -- The set of times for
observations $\mathscr{T}^{n} \equiv \mathscr{T}^{n}_{1}$ is a
renewal type process on $[0,+\infty[$ such that
\begin{equation*}
t_{0}=0 \quad \textrm{and} \quad t_{k}=\sum_{j=1}^{k} \tau_{j},
\quad 1 \leq k \leq n,
\end{equation*}
where $\tau^{n}=\{\tau_{k}, 1 \leq k \leq n\}$ is a sequence
of positive and i.i.d. random variables -- or \emph{inter-arrival times}
-- generated by a given probability density function $g(t)>0$ with
finite mean $\delta$. Let $g^{\star k}$ be the $k$th fold convolution of
$g$ with itself, then $g^{\star k}(t)$ is the density function of $t_{k}$ and
we define the \emph{renewal density} $h$ by $h(t):=\sum_{k=1}^{\infty}
g^{\star k}(t)$. Here and below, the function $h$ is supposed to be bounded by a constant $h_{0}$.

\begin{remark} \textnormal{The renewal density is known to satisfy
$h(u) \rightarrow \delta^{-1}$ as $u \rightarrow \infty$ (see Cox \cite{Cox1962}, p. 55), but its explicit expression is generally complicated to obtain. Nevertheless, the boundedness of $h$ is a condition which holds for a large class of
sequences $\tau^{n}$. For the reader convenience, we recall the example in Masry
\cite{Masry1983} corresponding to the usual situation where $\tau^{n}$ has a Gamma
density of type $r$, i.e.,
$$g(t) = \frac{(r/\delta) (r t /\delta)^{r-1} \exp(-r t/\delta)}{(r-1)!}, \quad r \in \mathbb{N} \backslash \{0\}, ~\delta > 0, ~t \geq 0,$$
with mean $\delta$ and variance $\delta^{2}/r$. Thus, if
$r=1$, $h(t)=\delta^{-1}$ for $t \geq 0$ ($\mathscr{T}^{n}_{1}$ is a Poisson process) and, if $r=2$, $h(t) = \delta^{-1} (1-\exp(-4 t/\delta))$ which approaches its limit $\delta^{-1}$ monotonically as $t \rightarrow \infty$. In both cases, the value
$h_{0}=\delta^{-1}$ is clearly appropriate. From $r=3$, the condition becomes delicate
to verify since $h(t)$ oscillates in approaching $\delta^{-1}$. The case $r=1$ is illustrated e.g. in A\"it-Sahalia and Mykland \cite{AitSahaliaMykland2003} with an example of financial data for which a histogram distribution of the sampling intervals is fitted by an exponential density.} \end{remark}

\noindent{\textbf{Jittered sampling}} -- First, we assume that
the process is regularly observed with a period $\delta > 0$. This
sequence $\mathscr{T}^{n} \equiv \mathscr{T}^{n}_{2}$ is then
contaminated by an additive noise to model the plausible
imperfections of a measurement recording system:
\begin{equation*}
t_{0}=Z_{0} \quad \textrm{and} \quad t_{k} = k \delta + Z_{k},
\quad 1 \leq k \leq n,
\end{equation*}
where $Z^{n} = \{Z_{k}, 0 \leq k \leq n\}$ denotes an i.i.d.
random sample from a symmetric probability density function
$g_{J}(z)$ over $[-\delta/2,\delta/2]$. In comparison with renewal times, jittered times could be seen as only partially random due to the
deterministic component in $t_{k}$.

\begin{remark} \textnormal{The observations drawn from each of these two random designs are by definition irregularly spaced in time, but the ``long-term'' expected inter-arrival time between two consecutive random instants is equal to $\delta$ in each case.} \end{remark}

Finally, we introduce a periodic scheme examined in Bosq
\cite{Bosq1997} for kernel density estimation where the sampling step $\delta_{n}$ is
$n$-decreasing in a deterministic manner. \vspace{.3cm}

\noindent{\textbf{High frequency sampling}} -- In order to
represent the occurrence of high frequency observations during a
long time, the sampling instants in $\mathscr{T}^{n} \equiv
\mathscr{T}^{n}_{3}$ are defined periodically as
\begin{equation*}
t_{0,n}=0 \quad \textrm{and} \quad t_{k,n} = k \delta_{n}, \quad 1
\leq k \leq n,
\end{equation*}
where $\delta_{n} > 0$ and $\delta_{n} \rightarrow 0^{+}$, $T_{n}
= n \delta_{n} \rightarrow \infty$ as $n \rightarrow \infty$. In the
sequel, we will give minimal thresholds $\delta_{n}^{*}$ over
which our estimators converge with the optimal rates of the i.i.d.
case. The knowledge of $\delta_{n}^{*}$ will also help us to
minimize the costs of estimation without altering the rates. To
explain, observe that two situations may occur in applications.
First, if the total time of observation is a given and large enough $T_{n}$, the value of a minimal $\delta_{n}^{*}$ allows to
select a maximal number $n^{*}$ of points in $[0,T_{n}]$ to
estimate $f$. On the other hand, consider that a maximal and large
enough sample size $n$ is available, then we can deduce from
$\delta_{n}^{*}$ a minimal sufficient time $T_{n}^{*}=n \delta_{n}^{*}$ of
observation (see Blanke and Pumo \cite{BlankePumo2003}).
Furthermore, we will emphasize the convenience of such a framework
to sample a continuous-time process. Thus, under the Castellana-Leadbetter's
conditions, i.e. $\int_{0}^{\infty} \sup_{x,y} |g_{u}(x,y)| \mathrm{d}u <
\infty$ and $g_{u}(\cdot,\cdot)$ is continuous at $(x,x)$ for any
$u>0$, Bosq \cite{Bosq1997} proved that $\delta_{n}$ can be chosen in order to obtain the full rate $1/T_{n}$ of the pointwise
mean squared error of kernel estimators. In that situation, the sampling scheme is said to be \emph{admissible}.
Concerning admissible sampling in nonparametric density
estimation, let us cite relevant works by Leblanc \cite{Leblanc1995}
for wavelets estimators, by Biau \cite{Biau2003} for spatial
kernel estimators, and by Comte and Merlevède
\cite{ComteMerlevede2005} and Blanke \cite{Blanke2006},
respectively, for projection and adaptive kernel estimators.

\subsection{Mean integrated squared error}
The global accuracy of density estimators can be measured by the
\emph{mean integrated squared error} which is the expected squared
distance between a density estimator $\hat{f}_{n}$ and the true
density $f$ integrated over $\mathbb{R}^{d}$:
\begin{equation*}
\mathrm{MISE}\big(\hat{f}_{n}\big) = \mathrm{E}
\int_{\mathbb{R}^{d}} \left(\hat{f}_{n}(x)-f(x)\right)^{2} \mathrm{d}x.
\end{equation*}
It is also the sum of the \emph{integrated squared bias} (ISB) and
the \emph{integrated variance} (IV):
\begin{equation*}
\mathrm{ISB}\big(\hat{f}_{n}\big) = \int_{\mathbb{R}^{d}}
\left(\mathrm{E}\big(\hat{f}_{n}(x)\big)-f(x)\right)^{2} \mathrm{d}x
\quad \textrm{and} \quad \mathrm{IV}\big(\hat{f}_{n}\big) =
\int_{\mathbb{R}^{d}} \mathrm{E}\left(\hat{f}_{n}(x)-\mathrm{E}\big(\hat{f}_{n}(x)\big)\right)^{2} \mathrm{d}x.
\end{equation*}

Let us fix the following usual notations: $\mathcal{C}^{k}\big(\mathbb{R}^{d}\big)$ denotes the
set of $k$-times continuously differentiable functions and
$L_{k}\big(\mathbb{R}^{d}\big)$ the set of functions with
integrable $k$th power over $\mathbb{R}^{d}$ such that
$\|f\|_{k}=(\int_{\mathbb{R}^{d}} |f(x)|^{k} \mathrm{d}x)^{1/k}$.

\section{Histogram} \label{histogram}
We primarily examine the \emph{histogram}, which is the oldest and
most popular nonparametric estimator. Because of its simplicity,
histogram is still widely used in presentation and
practice by statisticians. The theoretical properties have
been also extensively studied in the i.i.d. case and we may refer e.g. to Scott \cite{Scott1992} (Chapter 3) for
a background material. For continuous-time delivered observations, both optimal and full rates of MISE and
asymptotic normality under Castellana-Leadbetter's conditions are given in Lejeune \cite{Lejeune2006}. In this section we derive results for observations collected at discretized instants according to the sequences
$\mathscr{T}^{n}_{i}$, $i=1,2,3$.

\subsection{Definitions and assumptions} Prior to the definition of our
estimator, we introduce a partition of $\mathbb{R}^{d}$, say
$\Pi_{n}$, into hypercubes of volume $h_{n}^{d}$ such that $h_{n} \rightarrow 0^{+}$, $n h_{n}^{d}
\rightarrow \infty$ as $n \rightarrow \infty$:
$$\Pi_{n}=\big\{\pi_{nj},j \in \mathbb{Z}^{d}\big\},$$ and $$\pi_{nj}=\prod_{k=1}^{d}
\big[b_{j_{k}},b_{j_{k}+1}\big[ = \prod_{k=1}^{d}
\bigg[c_{j_{k}}-\frac{h_{n}}{2},c_{j_{k}}+\frac{h_{n}}{2}\bigg[,
\quad j=(j_{1},\ldots,j_{d})' \in \mathbb{Z}^{d},$$ where
$b_{j}=(b_{j_{1}},\ldots,b_{j_{d}})' \in \mathbb{R}^{d}$,
$b_{j_{k}+1}-b_{j_{k}}=h_{n}$ and $c_{j_{k}}=(b_{j_{k}}+b_{j_{k}+1})/2$. Here $h_{n}$ is the smoothing parameter commonly
referred to as the \emph{bin width}. Note that the extension to unequal bin sizes is straightforward with more notations. From now on, we will suppose for any $x \in \mathbb{R}^{d}$ the existence of an index $j(x,n)$
in $\mathbb{Z}^{d}$ such that $x \in \pi_{j(x,n)}$ ($=:\pi_{nj}$).

Given $\Pi_{n}$ and $n$ discretized observations $X_{t_{1}},\ldots,X_{t_{n}}$, the histogram estimator of $f$ is then defined as
\begin{equation*}
\hat{f}_{n}^{H}(x)=\sum_{j} \left[ \frac{1}{n h_{n}^{d}}
\sum_{k=1}^{n} \mathds{1}_{\pi_{nj}}(X_{t_{k}}) \right]
\mathds{1}_{\pi_{nj}}(x)=:\sum_{j} \hat{f}_{j}
\mathds{1}_{\pi_{nj}}(x), \quad x \in \mathbb{R}^{d},
\end{equation*}
where $\mathds{1}_{\pi_{nj}}$ denotes the indicator function of $\pi_{nj}$. In particular, $\hat{f}_{n}^{H}$ has a unique value,
denoted by $\hat{f}_{j}$, over each hypercube $\pi_{nj}$ of
$\Pi_{n}$, which explains its high computational advantage.

Let $\mathcal{A}$ and $\mathcal{B}$ be two sub-$\sigma$-algebras
of $\mathcal{F}$, we introduce the classical strong mixing
coefficient defined as
$$\alpha(\mathcal{A},\mathcal{B}):=\sup_{A \in \mathcal{A}, B \in
\mathcal{B}} |P(A \cap B)-P(A)P(B)|.$$
Let denote $\sigma(X)$ the $\sigma$-algebra of events generated by a random variable $X$. In the sequel, we will use the definition of a $2$-$\alpha$-mixing process $\{X_{t}, t \in \mathbb{R}\}$ given in Bosq \cite{Bosq1998} as
\begin{equation*}
\alpha_{X}^{(2)}(u):=\sup_{t \in \mathbb{R}} \alpha\big(\sigma(X_{t}),\sigma(X_{t+u})\big) \rightarrow 0 ~~\textrm{as} ~u \rightarrow \infty.
\end{equation*}
Note that such a condition only for the couples $(X_{t},X_{t+u})$ is less restrictive than the classical one introduced by Rosenblatt \cite{Rosenblatt1956}.

These are now the main assumptions over processes. \vspace{.3cm}

\noindent{\textbf{Assumptions $\mathbf{A_{0}}$}}
\begin{enumerate}
\item[(i)] $f \in \mathcal{C}^{2}\big(\mathbb{R}^{d}\big)$ so that all
the partial derivatives are square Riemann-integrable;
\item[(ii)] $f$ is continuous and $\|f\|_{\infty}=\sup_{y
\in \mathbb{R}^{d}} f(y) < \infty$.
\end{enumerate}

\newpage

\noindent{\textbf{Assumptions $\mathbf{A_{1}}$} (with renewal and jittered samplings)}
\begin{enumerate}
\item[(i)] There exists $u_{0}>0$ such that for any $u \geq u_{0}$: $\sup_{z \in \mathbb{R}^{d}} |g_{u}(y,z)| \leq
k(y)$ with $k(\cdot)$ a positive, continuous and integrable function
defined on $\mathbb{R}^{d}$; \item[(ii)] $X^{T}$
is an arithmetically strongly mixing (ASM) process i.e. there exists
$\rho > 2$, $a_{0}>0$ and $u_{1}>u_{0}$ such that for any $u \geq u_{1}$: $\alpha_{X}^{(2)}(u)=\alpha\big(\sigma(X_{0}),\sigma(X_{u})\big) \leq a_{0}
u^{-\rho}$.
\end{enumerate}

\noindent{\textbf{Assumptions $\mathbf{A'_{1}}$} (with high frequency sampling)}
\begin{enumerate}
\item[(i)] There exists $\gamma_{0}>0$ and $u_{0}>0$ such that for any $0< u \leq u_{0}$: $\forall y \in
\mathbb{R}^{d}$, $\sup_{z \in \mathbb{R}^{d}} f_{u}(y,z) \leq
\varphi(y) u^{-\gamma_{0}}$ with $\varphi(\cdot)$ a positive, continuous and integrable function
defined on $\mathbb{R}^{d}$;
\item[(ii)] There exists a positive, continuous and integrable function $k(\cdot)$
defined on $\mathbb{R}^{d}$ such that for any $u \geq u_{0}$: $ \forall y \in
\mathbb{R}^{d}$, $\sup_{z \in \mathbb{R}^{d}} |g_{u}(y,z)| \leq
k(y) \pi(u)$ where $\pi(\cdot)$ is a bounded and
ultimately decreasing function which satisfies
$\int_{u_{1}}^{\infty} \pi(u) \mathrm{d}u < \infty$, $u_{1}>u_{0}$.
\end{enumerate}

The assumptions above are classical in nonparametric estimation
with dependent data. $A_{0}$ displays some constraints of
regularity on the true density $f$. The condition $A_{0}(i)$ is specific to the bias
treatment, it was previously introduced by Lecoutre
\cite{Lecoutre1985} to study the multivariate histogram in the i.i.d. case.

The following conditions should take into account the local
behavior of sample paths as well as the properties of asymptotic
independence of processes (respectively described with the
behavior of $g_{u}$ for $u$ near the origin and for $u$ large).
$A_{1}(i)$ is a mild condition on $g_{u}$ for intermediate values of $u$. In particular, it
slightly weakens the assumption of boundedness on the conditional density used by Masry
\cite{Masry1983} and Carbon, Garel, and Tran
\cite{CarbonGarelTran1997}.

$A'_{1}(i)$ appears to be less usual in density estimation, but it
is a typical condition for the continuous-time framework to control the explosive behavior
of the joint densities $f_{u}(\cdot,\cdot)$ in a neighborhood of
$u=0$. Assumptions $A'_{1}$ are in the spirit of those made (and widely discussed) by Blanke and Pumo
\cite{BlankePumo2003}. Here $A'_{1}(i)$ is used with high frequency sampling to obtain optimal rates
together with a short sampling step $\delta_{n}$ depending on a positive known coefficient $\gamma_{0}$. Roughly speaking, the value of $\gamma_{0}$ is directly linked with the h\"olderian properties of sample paths and the dimension $d$: namely, one has $\gamma_{0}=d/2$ for a wide class of $d$-dimensional ergodic diffusion processes and $\gamma_{0}=d$ for ``smooth'' processes (see e.g. Blanke \cite{Blanke2006} for technical details).

Other assumptions, namely $A_{1}(ii)$ and $A'_{1}(ii)$,
ensure asymptotic independence between variables distant in time.
$A_{1}(ii)$ involves a mild version of $\alpha$-mixing which is well-known to be
weaker than many dependence structures as $\phi$, $\beta$ or
$\rho$-mixing (see e.g. Doukhan \cite{Doukhan1994}). Finally, admissible high frequency samplings are obtained under $A'_{1}(ii)$, a quite typical condition in this context.

\subsection{Rates of convergence}
Using each sampling design defined above, we will now establish the optimal rate of histograms. For the sake of
readability, some crucial lemmas which provide upper bounds for
the variances and the covariances of $\hat{f}_{n}^{H}$ are
postponed to the proofs. Let $f'_{i}:=\partial f/\partial x_{i}$
and define the \emph{roughness} $R$ of $f'_{i}$ by its squared
$L_{2}$-norm: $R\big(f'_{i}\big):=\int_{\mathbb{R}^{d}}
f'_{i}(x)^{2} \mathrm{d}x$. Since the bias of $\hat{f}_{n}^{H}$ only depends on the bin width and the true unknown density $f$,
and not on the dependence structure of the data, we recall the following result given by Lecoutre
\cite{Lecoutre1985} with multivariate independent observations.
\begin{lemma} \label{H-BCI} If Assumption $A_{0}(i)$ is satisfied then
\begin{equation*}
\mathrm{ISB}\big(\hat{f}_{n}^{H}\big) = \frac{h_{n}^{2}}{12}
R_{d}\big(f'\big) + o\big(h_{n}^{2}\big),
\end{equation*}
where $R_{d}\big(f'\big):=\sum_{i=1}^{d} R\big(f'_{i}\big)$.
\end{lemma}

\subsubsection{Renewal and jittered samplings}
Let us denote by $\lceil x \rceil$ the smallest integer not less than the real $x$. The first
part of the next theorem gives an asymptotic upper bound for IV.
Consequently, from an ad hoc choice of the bin width $h_{n}$ which
balances both ISB and IV terms, we infer that histograms can
reach the same optimal rate $n^{-2/(d+2)}$ of convergence to $f$
as in the i.i.d. case.

\begin{theorem} \label{H-MISE1}
\begin{enumerate}
\item Under $A_{0}(ii)$ and $A_{1}$ and if $f^{1-1/p} \in C^{1}\big(\mathbb{R}^{d}\big) \bigcap L_{1}\big(\mathbb{R}^{d}\big)$ for $1 < p < \rho-1$, then
\begin{equation*}
\limsup_{n \rightarrow \infty} n h_{n}^{d} ~\mathrm{IV}\big(\hat{f}_{n}^{H}\big) \leq 1 + C,
\end{equation*}
where $C=2 u_{0} h_{0}$ for the renewal sampling and $C=2 \left\lceil
\frac{u_{0}}{\delta} \right\rceil$ for the jittered sampling;
\item If in addition $A_{0}(i)$ holds then the choice $h_{n}=cn^{-1/(d+2)}$, $0<c<\infty$, yields
\begin{equation*}
\limsup_{n \rightarrow \infty} n^{\frac{2}{d+2}} ~\mathrm{MISE}\big(\hat{f}_{n}^{H}\big) \leq \frac{c^{2}}{12}
R_{d}\big(f'\big) + \frac{1}{c^{d}} (1 + C),
\end{equation*}
with same constant $C$.
\end{enumerate}
\end{theorem}

\begin{remark} \textnormal{If $p=\rho-1$ the rates in Theorem
\ref{H-MISE1} remain valid but with larger asymptotic constants
(see proofs). Thus, if for instance $\rho \geq 3$, one may choose
$p=2$ provided that $f^{1/2}$ is continuous and integrable.}
\end{remark}

\subsubsection{High frequency sampling}
The high frequency model is interesting to find some connections between both discrete and
continuous-time frameworks. Here the period $\delta_{n}$ is now a
function of the sampling size $n$ so that all observations can be
as close in time as desired provided $n$ large enough. Within this setup we also need to check the local condition $A'_{1}(i)$ on the joint density of $(X_{0},X_{u})$ for the small values of $u$, wherein a (known) coefficient $\gamma_{0}$ is linked with the regularity of sample paths. In this framework the previous optimal rate of order $n^{-2/(d+2)}$ is still preserved. Moreover, depending on the value $\gamma_{0}$, we can derive a minimal $\delta_{n}^{*}$ (more precisely $\delta_{n}^{*}(\gamma_{0})$) and then deduce a minimal time of observation of the process $T_{n}^{*}$ that ensures this rate.
\begin{theorem} \label{H-MISE2Opt} According to the value of $\gamma_{0}$ we assume that $\delta_{n} \geq \delta_{n}^{*}(\gamma_{0})$
defined as
\begin{equation} \label{DeltaOpt}
\delta_{n}^{*}(\gamma_{0}):=d_{1} h_{n}^{d}
\mathds{1}_{\{\gamma_{0}<1\}} + d_{2} h_{n}^{d}
\ln\big(h_{n}^{-d}\big) \mathds{1}_{\{\gamma_{0}=1\}} + d_{3}
h_{n}^{d/\gamma_{0}} \mathds{1}_{\{\gamma_{0}>1\}}, ~0<d_{1},
d_{2}, d_{3}<\infty. \end{equation}
\begin{enumerate} \item Then under $A_{0}(ii)$ and $A'_{1}$
\begin{equation*}
\limsup_{n \rightarrow \infty} n h_{n}^{d} ~\mathrm{IV}\big(\hat{f}_{n}^{H}\big) \leq 1+C_{\gamma_{0}},
\end{equation*}
where $C_{\gamma_{0}}$ is a positive constant which depends upon
$\gamma_{0}$ (see its explicit form in proofs);
\item If in addition $A_{0}(i)$ holds with $h_{n}=cn^{-1/(d+2)}$, $0<c<\infty$, then
\begin{equation*}
\limsup_{n \rightarrow \infty} n^{\frac{2}{d+2}} ~\mathrm{MISE}\big(\hat{f}_{n}^{H}\big) \leq \frac{c^{2}}{12}
R_{d}\big(f'\big) + \frac{1}{c^{d}}(1+C_{\gamma_{0}}),
\end{equation*}
with same constant $C_{\gamma_{0}}$.
\end{enumerate}
\end{theorem}

\begin{remark} \textnormal{Using $A_{0}$ with either $A_{1}$ or $A'_{1}$, our results in Theorems \ref{H-MISE1} and \ref{H-MISE2Opt} are similar to those derived with independent variables by Lecoutre \cite{Lecoutre1985} in the $d$-dimensional
setup. Thus we retrieve (in $\limsup$) the same optimal rate
$n^{-2/(d+2)}$ in terms of MISE. The additional asymptotic cons\-tant $C$ or $C_{\gamma_{0}}$ in the variance bound arises as a non negligible remainder of the covariance term; it clearly depends on the sampling design in use.
Nevertheless, if $\delta_{n}$ is such that $\delta_{n}/\delta_{n}^{*}(\gamma_{0}) \rightarrow \infty$ as $n
\rightarrow \infty$, we can remove $C_{\gamma_{0}}$ in
Theorem \ref{H-MISE2Opt} to get the exact limiting constant of the
i.i.d. case with $h_{n}=cn^{-1/(d+2)}$, $0<c<\infty$.}
\end{remark}

Remembering that $T_{n}=n \delta_{n}$ the rate $n^{-2/(d+2)}$ may be easily rewritten in terms of $T_{n}$ according to the value of $\gamma_{0}$.
\begin{corollary} \label{H-MISE2Full}
Under $A_{0}$ and $A'_{1}$ the choice $h_{n}=cn^{-1/(d+2)}$, $0<c<\infty$, leads to
\begin{equation*}
\mathrm{MISE}\big(\hat{f}_{n}^{H}\big) = \left\{
\begin{array}{ll}
O\big(T_{n}^{-1}\big) &  \textrm{with} ~\delta_{n}=d_{1} h_{n}^{d}, ~0 < d_{1} < \infty, ~\textrm{if} ~\gamma_{0}<1; \\
O\big(T_{n}^{-1} \ln T_{n}\big) & \textrm{with} ~\delta_{n}=d_{2} h_{n}^{d} \ln\big(h_{n}^{-d}\big), ~0 < d_{2} < \infty, ~\textrm{if} ~\gamma_{0}=1; \\ O\Big(T_{n}^{-\frac{2 \gamma_{0}}{2 \gamma_{0}+d(\gamma_{0}-1)}}\Big) & \textrm{with} ~\delta_{n}=d_{3} h_{n}^{d/\gamma_{0}}, ~0 < d_{3} < \infty,  ~\textrm{if} ~\gamma_{0}>1. \end{array} \right.
\end{equation*}
\end{corollary}

\begin{remark} \textnormal{For the special case of irregular paths processes ($\gamma_{0}<1$), we thus observe in Corollary \ref{H-MISE2Full} a surprising similarity between the best rate of order $1/T_{n}$ and the $1/T$-parametric rate encountered in the real continuous-time context. Indeed, the time of observation clearly depends on the value of $\gamma_{0}$ since $T_{n}$ has to be of order $n^{2/(d+2)}$ ($\gamma_{0}<1$), $n^{2/(d+2)} \ln n$ ($\gamma_{0}=1$) or $n^{(2 \gamma_{0}+d(\gamma_{0}-1))/((d+2) \gamma_{0})}$ ($\gamma_{0}>1$) so as to obtain same efficiency in estimation. Especially, this enlightens the fact that irregular paths processes may be observed less time than more regular ones ($\gamma_{0} \geq 1$).}
\end{remark}

Finally, it may be interesting to indicate the exact limit of the pointwise variance of $\hat{f}_{n}^{H}(x)$ in the case $\gamma_0<1$. The following proposition is thus obtained as a simple transposition from
kernel to histogram estimators of a result by Bosq \cite{Bosq1997}
(Proposition 7.1. (i)).
\begin{proposition} \label{H-Vexact} Let $x \in \mathbb{R}^{d}$ and assume
that
\begin{itemize}
\item[(i)] $\|g_{u}\|_{\infty} \leq \pi(u)$ where
$(1+u)\pi(u)$ is integrable over $]0,+\infty[$ and $u \pi(u)$
is bounded and ultimately decreasing. Furthermore
$g_{u}(\cdot,\cdot)$ is continuous at $(x,x)$;
\item[(ii)] $\sup_{(y,z) \in
\mathbb{R}^{2d}} \left|\sum_{r=1}^{\infty} \delta_{n} g_{r
\delta_{n}}(y,z)-\int_{0}^{\infty} g_{u}(y,z) \mathrm{d}u\right|
\rightarrow 0$ as $\delta_{n} \downarrow 0^{+}$,
\end{itemize} then
\begin{equation*} \lim_{n \rightarrow \infty} T_{n} ~\mathrm{Var}\big(\hat{f}_{n}^{H}(x)\big) = 2 \int_{0}^{\infty}
g_{u}(x,x) \mathrm{d}u,
\end{equation*}
provided that $\delta_{n}=o\big(h_{n}^{d}\big)$.
\end{proposition}
\begin{remark} \label{remark6} \textnormal{From Kutoyants \cite{Kutoyants1998}, the limiting constant is also the minimax bound for mean squared error in the case of ergodic diffusion processes satisfying some conditions of regularity on the trend coefficient and the diffusion coefficient (see Veretennikov \cite{Veretennikov1999}).}
\end{remark}

\section{Frequency polygon} \label{FreqPoly}
Given a (univariate) histogram, the \emph{frequency polygon} results from a
natural smoothing with straight lines to get a continuous estimator. However, the gain of this simple linear smoothing is substantial since we
immediately improve the weak order $h_{n}^{2}$ inherent to the
bias of histograms. The main properties of frequency polygons are gathered in Scott
\cite{Scott1992} (Chapter 4) within the i.i.d. setup. The
mixing case was then treated by Carbon, Garel, and Tran
\cite{CarbonGarelTran1997}, and recently extended to the random
fields by Carbon \cite{Carbon2008}. In continuous-time, Lejeune
\cite{Lejeune2005,Lejeune2006} established both optimal and parametric
rates of MISE and asymptotic normality; the extension to the
random fields is done in a submitted work by Bensa\"id and
Dabo-Niang \cite{BensaidDabo2007}. For the sake of simplicity, we shall confine attention to
the real case ($d=1$ and $\gamma_{0} \leq 1$).

For convenience, $f'$ and $f''$ denote the first and second derivatives of $f$ and we define the roughness of $f''$ by $R\big(f''\big):=\int_{\mathbb{R}} f''(x)^{2} \mathrm{d}x$.

\subsection{Definition and assumptions} Based upon $\Pi_{n}$ and $X_{t_{1}},\ldots,X_{t_{n}}$, the frequency polygon is simply constructed by connecting the mid-points of the histogram heights with straight line segments
\begin{equation*}
\hat{f}_{n}^{FP}(x) = \sum_{j} \left[ \left(
\frac{x-c_{j}}{h_{n}} \right) \hat{f}_{j+1} + \left(
\frac{c_{j+1}-x}{h_{n}} \right) \hat{f}_{j} \right]
\mathds{1}_{[c_{j},c_{j+1}[}(x), \quad x \in \mathbb{R}.
\end{equation*}

In the literature we find also alternative definitions which differ
from the way of interpolation as e.g. the \emph{edge frequency
polygon} introduced by Jones, Samiuddin and Al-Harbey Maatouk
\cite{Jones1998} or its extended form by Dong and Zheng
\cite{DongZheng2001}. All these estimators share the same rates of convergence but with different asymptotic constants.

In agreement with assumptions $A_{1}$ and $A'_{1}$ in the previous section we will describe the proper\-ties of the frequency polygon under the following conditions on $f$. \vspace{.3cm}

\noindent{\textbf{Assumptions $\mathbf{A'_{0}}$}}
\begin{enumerate}
\item[(i)] $f \in \mathcal{C}^{2}(\mathbb{R})$, $f'' \in L_{1}(\mathbb{R})$ and $f, f'' \in
L_{2}(\mathbb{R})$; \item[(ii)] $\left|f''(x)-f''(y)\right| \leq
l_{0} |x-y|^{\nu}$, $l_{0}>0$, $\nu \in ]0,1]$, for $(x,y) \in
\mathbb{R}^{2}$;
\item[(iii)] $f$ is continuous and $\|f\|_{\infty} < \infty$.
\end{enumerate}

\subsection{Rates of convergence} The ISB
contribution is given in Scott \cite{Scott1985}.
\begin{lemma} \label{FP-BCI}
If Assumptions $A'_{0}(i)(ii)$ are satisfied then
\begin{equation*}
\mathrm{ISB}\big(\hat{f}_{n}^{FP}\big) = \frac{49}{2880}
R\big(f''\big) h_{n}^{4} + o\big(h_{n}^{4}\big).
\end{equation*}
\end{lemma}

\begin{remark} \textnormal{The nice order $h_{n}^{4}$ is much better compared with histograms and familiar for more sophisticated density estimators as kernel estimators. As emphasized earlier the bias term does not depend on the sampling scheme.}
\end{remark}

\subsubsection{Renewal and jittered samplings}
Using the analysis on histograms with a new suitable choice of $h_{n}$ we give the
optimal rate of frequency polygons. Note that constants $C$ and
$C_{\gamma_{0}}$ are unchanged.
\begin{theorem} \label{FP-MISE1}
\begin{enumerate}
\item Under $A'_{0}(iii)$ and $A_{1}$ and if $f^{1-1/p} \in C^{1}(\mathbb{R}) \bigcap L_{1}(\mathbb{R})$ for $1 < p < \rho-1$, then
\begin{equation*}
\limsup_{n \rightarrow \infty} n h_{n} ~\mathrm{IV}\big(\hat{f}_{n}^{FP}\big) \leq \frac{2}{3} + C;
\end{equation*}
\item If in addition $A'_{0}(i)(ii)$ hold then the choice $h_{n}=cn^{-1/5}$, $0<c<\infty$, yields
\begin{equation*}
\limsup_{n \rightarrow \infty} n^{\frac{4}{5}} ~\mathrm{MISE}\big(\hat{f}_{n}^{FP}\big) \leq \frac{49}{2880} c^{4}
R\big(f''\big) + \frac{1}{c} \left( \frac{2}{3} + C \right).
\end{equation*}
\end{enumerate}
\end{theorem}

\subsubsection{High frequency sampling}
Finally, recovering the local properties of sample paths when data become dense in time, we
find again the optimal rate while minimizing the time of observation.
\begin{theorem} \label{FP-MISE2Opt} According to the values of $\gamma_{0}$, we consider optimal choices $\delta_{n}^{*}(\gamma_{0})$
given by (\ref{DeltaOpt}).
\begin{enumerate}
\item Then under $A'_{0}(iii)$ and $A'_{1}$
\begin{equation*}
\limsup_{n \rightarrow \infty} n h_{n} ~\mathrm{IV}\big(\hat{f}_{n}^{FP}\big) \leq
\frac{2}{3}+C_{\gamma_{0}};
\end{equation*}
\item If in addition $A'_{0}(i)(ii)$ hold with $h_{n}=cn^{-1/5}$, $0 < c < \infty$, then
\begin{equation*}
\limsup_{n \rightarrow \infty} n^{\frac{4}{5}} ~\mathrm{MISE}\big(\hat{f}_{n}^{FP}\big) \leq \frac{49}{2880} c^{4}
R\big(f''\big) + \frac{1}{c} \left(\frac{2}{3} +
C_{\gamma_{0}}\right).
\end{equation*}
\end{enumerate}
\end{theorem}

\begin{remark} \textnormal{In both Theorems \ref{FP-MISE1} and \ref{FP-MISE2Opt} we exhibit (in $\limsup$) the same
$n^{-4/5}$-consistency obtained in Scott \cite{Scott1985} with i.i.d. observations.
The additional asymptotic constant $C$ or $C_{\gamma_{0}}$ still stays and relies on the sampling design in use; but, in Theorem \ref{FP-MISE2Opt}, any choice of $\delta_{n}$ satisfying $\delta_{n}/\delta_{n}^{*}(\gamma_{0}) \rightarrow \infty$ as $n \rightarrow \infty$ allows to remove $C_{\gamma_{0}}$ to get the exact limiting constant of the i.i.d. case with $h_{n}=c n^{-1/5}$. Finally, note that if we take $p=\rho-1$ in Theorem \ref{FP-MISE1} the rates remain valid up to a decayed constant.}
\end{remark}

\begin{corollary} \label{FP-MISE2Full}
Under $A'_{0}$ and $A'_{1}$ the choice $h_{n}=cn^{-1/5}$, $0<c<\infty$,  leads to
\begin{equation*}
\mathrm{MISE}\big(\hat{f}_{n}^{FP}\big) = \left\{
\begin{array}{ll}
O\big(T_{n}^{-1}\big) &  \textrm{with} ~\delta_{n} = d_{1} h_{n}, ~0 < d_{1} < \infty, ~\textrm{if} ~\gamma_{0}<1; \\
O\big(T_{n}^{-1} \ln T_{n}\big) & \textrm{with} ~\delta_{n} = d_{1} h_{n} \ln\big(h_{n}^{-1}\big), ~0 < d_{2} < \infty, ~\textrm{if} ~\gamma_{0}=1.
\end{array} \right.
\end{equation*}
\end{corollary}

\begin{remark} \textnormal{As noticed before real irregular paths processes may be observed less time than more regular ones since $T_{n}$ has to be of order $n^{4/5}$ ($\gamma_{0}<1$) or $n^{4/5} \ln n$ ($\gamma_{0}=1$) to obtain same efficiency in estimation.}
\end{remark}

For completeness, the exact limit of the pointwise variance of
$\hat{f}_{n}^{FP}(x)$ follows straightforwardly from Proposition
\ref{H-Vexact} (see also Remark \ref{remark6}).
\begin{proposition} \label{FP-Vexact} Under conditions of Proposition \ref{H-Vexact} with $\delta_{n}=o(h_{n})$, one has
\begin{equation*} \lim_{n \rightarrow \infty} T_{n} ~\mathrm{Var}\big(\hat{f}_{n}^{FP}(x)\big) = 2 \int_{0}^{\infty}
g_{u}(x,x) \mathrm{d}u, \quad x \in \mathbb{R}.
\end{equation*}
\end{proposition}

\section{Discussion} \label{conclusion}
In this work we derive the optimal $L_{2}$-rates of two
computationally advantageous density estimators in the setup where observations are discretely sampled from a
continuous-time process. For practical considerations we have studied three time sampling procedures to properly describe the time occurrences of the real data. Thus, values may be available at low or high frequency but also regularly or irregularly spaced in time. Therefore our main results state that all designs either random or deterministic lead to the optimal rates $n^{-2/(d+2)}$ for histograms and $n^{-4/5}$ ($d=1$) for frequency polygons, with respect to the MISE convergence, which are those derived in the i.i.d. case. From this result, the frequency polygon is a good alternative to more sophisticated nonparametric density estimators. Particularly, we have focused on a high frequency sampling to reveal some parallels with the idealized continuous-time framework as soon as observations are selected close enough to each other. We then use the local properties of sample paths to have a consistent estimation with a minimal time of experiment. This fact might be explained as follows:
irregular sample paths carry much more information than regular ones where the correlation between two successive variables is
much stronger. Consequently, we infer that more the paths are irregular -- i.e. when $A'_{1}(i)$ holds with $\gamma_{0}<1$ -- more the time of observation would be shortened with a good behavior of the both estimators.
Although not presented here, simulations in progress already corroborate our theoretical results in the particular case of two stationary real gaussian processes. As awaited the frequency polygon performs well and appears much closer to kernel estimator than to histogram. To go further in our investigations, it remains to examine the case of non-gaussian processes including, for instance, the cumbersome problem of estimating bimodal densities. The important issue of finding optimal choices for the bin width value is left for future work.

\section{Proofs} \label{Proofs}
Throughout this section, we detail the proofs of Theorems
\ref{H-MISE1}, \ref{H-MISE2Opt}, \ref{FP-MISE1} and
\ref{FP-MISE2Opt}. In order to do this, some auxiliary lemmas are
necessary to derive upper bound expressions for the variance of
$\hat{f}_{n}^{H}(x)$, $x \in \pi_{nj}$, which will depend on the
sampling scheme being used. Let $\|X\|_{q} =
(\mathrm{E}|X|^{q})^{1/q}$ with $1 \leq q < \infty$, then $X \in
L_{q}(P)$ means that $\|X\|_{q} < \infty$. We recall the following
useful covariance inequality as written in Bosq \cite{Bosq1998}
(p. 21).
\begin{lemma}[Davydov's inequality] \label{Davydov}
Let $X \in L_{q}(P)$ and $Y \in L_{r}(P)$ with $q>1$, $r>1$ and $\frac{1}{q}+\frac{1}{r}<1$, then
\begin{equation*}
|\mathrm{Cov}(X,Y)| \leq 2 p \Big[2 \alpha
\big(\sigma(X),\sigma(Y)\big)\Big]^{1/p} \|X\|_{q} \|Y\|_{r},
\end{equation*}
where $\frac{1}{p}+\frac{1}{q}+\frac{1}{r}=1$.
\end{lemma}

\subsection{Histogram}

\subsubsection{Variance bounds with random sampling}

\begin{lemma}[\textnormal{\bf renewal sampling}] \label{H-MajVar1} If $A_{0}(ii)$ and $A_{1}$ hold then we obtain for $1 < p \leq \rho-1$:
\begin{multline} \label{VarBound1} n h_{n}^{d} ~\mathrm{Var}\big(\hat{f}_{j}\big) \leq f(\xi_{j}) (1 - h_{n}^{d}
f(\xi_{j})) (1+2 u_{0} h_{0}) + 2 h_{0} k\big(\dot{\xi}_{j}\big) h_{n}^{\varepsilon} \\
+ \frac{4 p^{2} (2 a_{0})^{1/p} h_{0}}{\rho - p}
f(\xi_{j})^{1-\frac{1}{p}} ~h_{n}^{\frac{1}{p}
\{(d-\varepsilon)(\rho-p)-d \}},
\end{multline}
with $0 \leq \varepsilon \leq d \left(1 - \frac{1}{\rho - p}
\right)$ and $\big(\xi_{j},\dot{\xi}_{j}\big) \in \pi_{nj}^{2}$.
\end{lemma}

\begin{lemma}[\textnormal{\bf jittered sampling}] \label{H-MajVar2} Under the same conditions as in Lemma \ref{H-MajVar1} and $1 < p \leq \rho-1$:
\begin{multline} \label{VarBound2}
n h_{n}^{d} ~\mathrm{Var}\big(\hat{f}_{j}\big) \leq
f(\xi_{j}) (1-h_{n}^{d} f(\xi_{j})) \left(1 + 2 \left\lceil
\frac{u_{0}}{\delta} \right\rceil \right) + 2 k\big(\dot{\xi}_{j}\big) \left(h_{n}^{\varepsilon}-\left\lceil
\frac{u_{0}}{\delta} \right\rceil h_{n}^{d}\right) \\
+ \frac{4 p^{2} (2 a_{0})^{1/p}}{(\rho - p)
\delta^{\frac{\rho}{p}}} f(\xi_{j})^{1-\frac{1}{p}} ~h_{n}^{\frac{1}{p} \{(d-\varepsilon)(\rho-p)-d \}} \big(1-2
h_{n}^{d-\varepsilon}\big)^{1-\frac{\rho}{p}},
\end{multline}
with $0 \leq \varepsilon \leq d \left( 1 - \frac{1}{\rho-p}
\right)$ and $\big(\xi_{j},\dot{\xi}_{j}\big) \in \pi_{nj}^{2}$.
\end{lemma}

For further use, we give the proofs for the covariances. \vspace{.3cm}

\noindent \textbf{Proof of Lemma \ref{H-MajVar1}} ~For any $(x,y)
\in \mathbb{R}^{d} \times \mathbb{R}^{d}$, we suppose the
existence of two indexes $j_{1}(x,n)$ and $j_{2}(y,n)$ in $\mathbb{Z}^{d}$
such that $x \in \pi_{j_{1}(x,n)} (=:\pi_{nj_{1}})$ and $y \in
\pi_{j_{2}(y,n)} (=:\pi_{nj_{2}})$. Thus
\begin{equation*} \hat{f}_{n}^{H}(x)=\hat{f}_{j_{1}}=\frac{1}{n h_{n}^{d}}
\sum_{k=1}^{n} \mathds{1}_{\pi_{nj_{1}}}(X_{t_{k}}), \quad
\hat{f}_{n}^{H}(y)=\hat{f}_{j_{2}}=\frac{1}{n h_{n}^{d}}
\sum_{k=1}^{n} \mathds{1}_{\pi_{nj_{2}}}(X_{t_{k}}), \end{equation*} and
\begin{align*} n h_{n}^{d} ~\mathrm{Cov}\big(\hat{f}_{j_{1}},\hat{f}_{j_{2}}\big) = &
\frac{1}{n h_{n}^{d}} \sum_{k=1}^{n}
\mathrm{Cov}\big(\mathds{1}_{\pi_{nj_{1}}}(X_{t_{k}}),\mathds{1}_{\pi_{nj_{2}}}(X_{t_{k}})
\big) \\ & + \frac{2}{n h_{n}^{d}} \sum_{p=1}^{n-1}
\sum_{q=p+1}^{n}
\mathrm{Cov}\big(\mathds{1}_{\pi_{nj_{1}}}(X_{t_{p}}),\mathds{1}_{\pi_{nj_{2}}}(X_{t_{q}})\big)
=: V_{n} + C_{n}.
\end{align*}
Set $p_{k}:=P(X_{0} \in \pi_{nk})$, $k \in \mathbb{Z}^{d}$. The
``variance term'' $V_{n}$ is easy to compute.
\begin{equation*}
V_{n} = \frac{1}{n h_{n}^{d}} \sum_{k=1}^{n}
\mathrm{Cov}\big(\mathds{1}_{\pi_{nj_{1}}}(X_{0}),\mathds{1}_{\pi_{nj_{2}}}(X_{0})\big)
= \frac{1}{h_{n}^{d}} \big(P(X_{0} \in \pi_{nj_{1}}, X_{0} \in
\pi_{nj_{2}}) - p_{j_{1}} p_{j_{2}}\big).
\end{equation*}
Since $f$ is continuous there exists at least one point $\xi_{j}
\in \pi_{nj}$ such that $\int_{\pi_{nj}} f(x) \mathrm{d}x = h_{n}^{d}
f(\xi_{j})$. Then if $j_{1} \neq j_{2}$, we get
$$V_{n} = -\frac{1}{h_{n}^{d}} p_{j_{1}} p_{j_{2}} = -h_{n}^{d} f(\xi_{j_{1}}) f(\xi_{j_{2}}),$$
where $\big(\xi_{j_{1}},\xi_{j_{2}}\big) \in \pi_{nj_{1}} \times
\pi_{nj_{2}}$. Otherwise if $j_{1}=j_{2}=j$:
$$V_{n} = \frac{1}{h_{n}^{d}} p_{j} (1 - p_{j}) = f(\xi_{j}) (1 - h_{n}^{d} f(\xi_{j})).$$
Let us turn to the ``covariance term'' $C_{n}$. By stationarity and, since $t_{p}-t_{q}$ and $t_{p-q}$ are equal in distribution, we
have
\begin{align*} C_{n} & = \frac{2}{n h_{n}^{d}} \sum_{r=1}^{n-1} \sum_{p=1}^{n-r}
\mathrm{Cov}\big(\mathds{1}_{\pi_{nj_{1}}}(X_{0}),\mathds{1}_{\pi_{nj_{2}}}(X_{t_{p+r}-t_{p}})\big) \\
& = \frac{2}{h_{n}^{d}} \sum_{r=1}^{n-1} \left(1-\frac{r}{n}\right) \int_{0}^{\infty} \mathrm{Cov}\big(\mathds{1}_{\pi_{nj_{1}}}(X_{0}),\mathds{1}_{\pi_{nj_{2}}}(X_{u})\big) g^{\star r}(u) \mathrm{d}u =: C_{n,1} + C_{n,2} + C_{n,3},
\end{align*}
where \begin{equation*} C_{n,i} := \frac{2}{h_{n}^{d}}
\sum_{r=1}^{n-1} \left(1-\frac{r}{n}\right) \int_{E_{i}}
\mathrm{Cov}\big(\mathds{1}_{\pi_{nj_{1}}}(X_{0}),\mathds{1}_{\pi_{nj_{2}}}(X_{u})
\big) g^{\star r}(u) \mathrm{d}u, ~~i=1,2,3, \end{equation*} with $E_{1}=(0,u_{0})$, $E_{2}=\big(u_{0},h_{n}^{-d+\varepsilon}\big)$ and $E_{3}=\big(h_{n}^{-d+\varepsilon},\infty\big)$, for some $0 \leq \varepsilon < d$
to be specified later. Recall that
$h(u)=\sum_{r=1}^{\infty} g^{\star r}(u)$, one seeks to bound each
covariance subterm. First, by Cauchy-Schwarz inequality and Fubini's
theorem,
\begin{align*} |C_{n,1}| & \leq \frac{2}{h_{n}^{d}} \sqrt{\mathrm{Var}\big(\mathds{1}_{\pi_{nj_{1}}}(X_{0})\big)}
\sqrt{\mathrm{Var}\big(\mathds{1}_{\pi_{nj_{2}}}(X_{0})\big)} \int_{0}^{u_{0}} h(u) \mathrm{d}u \\
& \leq 2 u_{0} h_{0} \sqrt{f(\xi_{j_{1}})
f(\xi_{j_{2}}) (1-h_{n}^{d} f(\xi_{j_{1}})) (1-h_{n}^{d}
f(\xi_{j_{2}}))}.
\end{align*}
Then $A_{1}(i)$ and Fubini imply
\begin{equation*}
|C_{n,2}| \leq \frac{2}{h_{n}^{d}} \int_{u_{0}}^{h_{n}^{-d+\varepsilon}} \iint_{\pi_{nj_{1}} \times \pi_{nj_{2}}} \sup_{y \in \mathbb{R}^{d}} |g_{u}(x,y)| \mathrm{d}x \mathrm{d}y ~h(u) \mathrm{d}u \leq 2 h_{0} k\big(\dot{\xi}_{j_{1}}\big) h_{n}^{\varepsilon},
\end{equation*}
where $\dot{\xi}_{j_{1}} \in \pi_{nj_{1}}$. \newline

\noindent Now, it is clear that for $n$ large enough we have
$h_{n}^{-d+\varepsilon} \geq u_{1}$. So using Davydov's inequality
(Lemma \ref{Davydov}) with mixing condition $A_{1}(ii)$ and
Fubini, for any $(p,q) \in \big]1,\rho-1\big] \times
\left[1+\frac{\rho}{\rho-2},\infty\right[$ such that $\frac{2}{q}
+ \frac{1}{p}=1$:
\begin{align*}
|C_{n,3}| & \leq \frac{2}{h_{n}^{d}} \int_{h_{n}^{-d +
\varepsilon}}^{\infty} 2p ~2^{1/p} \|
\mathds{1}_{\pi_{nj_{1}}}(X_{0}) \|_{q} \|
\mathds{1}_{\pi_{nj_{2}}}(X_{u}) \|_{q} \big(\alpha_{X}^{(2)}(u)\big)^{1/p} ~h(u) \mathrm{d}u \\
& \leq 4p (2 a_{0})^{1/p} h_{0}
\sqrt{f(\xi_{j_{1}})^{1-\frac{1}{p}} f(\xi_{j_{2}})^{1-\frac{1}{p}}} ~h_{n}^{-\frac{d}{p}}
\int_{h_{n}^{-d+\varepsilon}}^{\infty} u^{-\frac{\rho}{p}} \mathrm{d}u \\
& \leq \frac{4 p^{2} (2 a_{0})^{1/p} h_{0}}{\rho-p}
\sqrt{f(\xi_{j_{1}})^{1-\frac{1}{p}} f(\xi_{j_{2}})^{1-\frac{1}{p}}} ~h_{n}^{\frac{1}{p} \{(d-\varepsilon)(\rho-p)-d\}}.
\end{align*}
Finally, setting $k_{1}:=2 u_{0} h_{0}$, $k_{2}:=2 h_{0}$ and $k_{3}:=\frac{4 p^{2} (2
a_{0})^{1/p} h_{0}}{\rho-p}$, one has
\begin{multline} \label{CovBound1}
n h_{n}^{d} ~\mathrm{Cov}\big(\hat{f}_{j_{1}},\hat{f}_{j_{2}}\big) \leq
-h_{n}^{d} f(\xi_{j_{1}}) f(\xi_{j_{2}}) + k_{1}
\sqrt{f(\xi_{j_{1}}) f(\xi_{j_{2}}) (1-h_{n}^{d} f(\xi_{j_{1}})) (1-h_{n}^{d} f(\xi_{j_{2}}))} \\
+ k_{2} k\big(\dot{\xi}_{j_{1}}\big) h_{n}^{\varepsilon} +
k_{3} \sqrt{f(\xi_{j_{1}})^{1-\frac{1}{p}}
f(\xi_{j_{2}})^{1-\frac{1}{p}}} ~h_{n}^{\frac{1}{p}
\{(d-\varepsilon)(\rho-p)-d\}}.
\end{multline}
We then deduce the lemma by taking $j_{1}=j_{2}=j$ with the appropriate
expression of $V_{n}$. It turns out that the covariance is a
$O\big(1/\big(n h_{n}^{d}\big)\big)$ for any choice of
$\varepsilon$ pertaining to $\left[0, d \Big(1-\frac{1}{\rho-p}
\Big)\right]$. ~~$\blacksquare$ \newline

\noindent \textbf{Proof of Lemma \ref{H-MajVar2}} ~Here the
calculus of $V_{n}$ is exactly the same as in the proof of Lemma
\ref{H-MajVar1}. In fact, the delicate point will consist again in
bounding $C_{n}$. To do so, we give the common probability density
function, say $\Delta_{Z}$, of all random variables $\{Z_{j} -
Z_{i}$, $i < j\}$. Since the variables $\{Z_{i}$, $0 \leq i \leq
n\}$ are supposed to be independent and symmetrically distributed,
we have $\Delta_{Z}(t) = {g_{J}}^{\star 2}(t) = \int_{\mathbb{R}}
g_{J}(t-y) g_{J}(y) \mathrm{d}y$ with support over $[-\delta,\delta]$. Let
us denote by $\lfloor x \rfloor$ the largest integer less than or
equal to the real $x$, and set $r_{0}:=\lceil u_{0} / \delta \rceil$
and $r_{n}^{1}:=\lfloor h_{n}^{-d+\varepsilon} \rfloor$ for some
$0 \leq \varepsilon < d$ to be specified later. Now stationarity implies
\begin{align*} C_{n} & = \frac{2}{n h_{n}^{d}} \sum_{r=1}^{n-1} \sum_{p=1}^{n-r}
\mathrm{Cov}\big(\mathds{1}_{\pi_{nj_{1}}}(X_{0}),\mathds{1}_{\pi_{nj_{2}}}(X_{t_{p+r}-t_{p}})\big) \\
& = \frac{2}{h_{n}^{d}} \sum_{r=1}^{n-1} \left(1 -
\frac{r}{n}\right) \int_{-\delta}^{\delta} \mathrm{Cov}\big(\mathds{1}_{\pi_{nj_{1}}}(X_{0}),\mathds{1}_{\pi_{nj_{2}}}(X_{r \delta + t})\big) \Delta_{Z}(t) \mathrm{d}t =: C_{n,1} + C_{n,2} + C_{n,3},
\end{align*}
where \begin{equation*} C_{n,i} := \frac{2}{h_{n}^{d}}
\sum_{r \in E_{i}} \left(1-\frac{r}{n}\right)
\int_{-\delta}^{\delta}
\mathrm{Cov}\big(\mathds{1}_{\pi_{nj_{1}}}(X_{0}),\mathds{1}_{\pi_{nj_{2}}}(X_{r
\delta + t})\big) \Delta_{Z}(t) \mathrm{d}t, ~~i=1,2,3, \end{equation*} with $E_{1}=\{1,\ldots,r_{0}\}$, $E_{2}=\big\{r_{0}+1,\ldots,r_{n}^{1}\big\}$ and $E_{3}=\big\{r_{n}^{1}+1,\ldots,n-1\big\}$. By Cauchy-Schwarz
we get
\begin{align*}
|C_{n,1}| & \leq \frac{2}{h_{n}^{d}} \sum_{r=1}^{r_{0}}
\sqrt{\mathrm{Var}\big(\mathds{1}_{\pi_{nj_{1}}}(X_{0})\big)}
\sqrt{\mathrm{Var}\big(\mathds{1}_{\pi_{nj_{2}}}(X_{0})\big)} \int_{-\delta}^{\delta} \Delta_{Z}(t) \mathrm{d}t \\
& \leq 2 \left\lceil \frac{u_{0}}{\delta} \right\rceil
\sqrt{f(\xi_{j_{1}}) f(\xi_{j_{2}}) (1-h_{n}^{d} f(\xi_{j_{1}}))
(1-h_{n}^{d} f(\xi_{j_{2}}))}.
\end{align*}
Then using $A_{1}(i)$
\begin{align*}
|C_{n,2}| & \leq \frac{2}{h_{n}^{d}} \sum_{r=r_{0}+1}^{r_{n}^{1}}
\int_{-\delta}^{\delta} \Delta_{Z}(t)
\iint_{\pi_{nj_{1}} \times \pi_{nj_{2}}} \sup_{y \in \mathbb{R}^{d}} |g_{r \delta + t}(x,y)| \mathrm{d}x \mathrm{d}y ~\mathrm{d}t \\
& \leq 2 (r_{n}^{1}-r_{0}) h_{n}^{d} k\big(\dot{\xi}_{j_{1}}\big) \int_{-\delta}^{\delta} \Delta_{Z}(t) \mathrm{d}t \\
& \leq 2 k\big(\dot{\xi}_{j_{1}}\big) \left(
h_{n}^{-d+\varepsilon} - \left\lceil \frac{u_{0}}{\delta}
\right\rceil \right) h_{n}^{d}.
\end{align*}
By Davydov's inequality and $A_{1}(ii)$,
$$|C_{n,3}| \leq \frac{2}{h_{n}^{d}} \sum_{r=r_{n}^{1}+1}^{n-1} \int_{-\delta}^{\delta}
\left|\mathrm{Cov}\big(\mathds{1}_{\pi_{nj_{1}}}(X_{0}),\mathds{1}_{\pi_{nj_{2}}}(X_{r
\delta + t})\big)\right| \Delta_{Z}(t) \mathrm{d}t.$$ For any $(p,q) \in
\big]1,\rho-1 \big] \times \left[1 + \frac{\rho}{\rho-2},\infty
\right[$ such that $\frac{2}{q} + \frac{1}{p}=1$ and since
$\alpha_{X}^{(2)}(\cdot)$ is arithmetically decreasing, we have
\begin{align*}
|C_{n,3}| & \leq \frac{1}{h_{n}^{d}} ~4p ~2^{1/p} ~h_{n}^{\frac{2d}{q}} f(\xi_{j_{1}})^{1/q}
f(\xi_{j_{2}})^{1/q} \sum_{r=r_{n}^{1}+1}^{n-1} \int_{-\delta}^{\delta} \Big(\alpha_{X}^{(2)}\big((r-1) \delta\big)\Big)^{1/p} \Delta_{Z}(t) \mathrm{d}t \\
& \leq h_{n}^{-\frac{d}{p}} ~\frac{4 p ~2^{1/p}}{\delta}
\sqrt{f(\xi_{j_{1}})^{1-\frac{1}{p}} f(\xi_{j_{2}})^{1-\frac{1}{p}}} ~\int_{(r_{n}^{1}-1) \delta}^{\infty}
\big(\alpha_{X}^{(2)}(u)\big)^{1/p} \mathrm{d}u \\
& \leq \frac{4 p^{2} (2 a_{0})^{1/p}}{(\rho-p)
\delta^{\frac{\rho}{p}}} \sqrt{f(\xi_{j_{1}})^{1-\frac{1}{p}}
f(\xi_{j_{2}})^{1-\frac{1}{p}}} ~h_{n}^{-\frac{d}{p}}
\big(r_{n}^{1}-1\big)^{1-\frac{\rho}{p}}.
\end{align*}
Now if $p<\rho$ $\big(1-\frac{\rho}{p}<0\big)$ and since
$r_{n}^{1}>h_{n}^{-d+\varepsilon}-1$ we may write
\begin{equation*}
|C_{n,3}| \leq \frac{4 p^{2} (2 a_{0})^{1/p}}{(\rho-p)
\delta^{\frac{\rho}{p}}} \sqrt{f(\xi_{j_{1}})^{1-\frac{1}{p}}
f(\xi_{j_{2}})^{1-\frac{1}{p}}} ~h_{n}^{\frac{1}{p}
\{(d-\varepsilon)(\rho-p)-d\}} \big(1-2 h_{n}^{d-\varepsilon}\big)^{1-\frac{\rho}{p}},
\end{equation*}
where $\left(1-2 h_{n}^{d-\varepsilon}\right)^{1-\frac{\rho}{p}}
\rightarrow 1$ as $n \rightarrow \infty$. Hence we obtain
\begin{equation*}
|C_{n,3}| \leq \frac{4 p^{2} (2 a_{0})^{1/p}}{(\rho-p)
\delta^{\frac{\rho}{p}}} \sqrt{f(\xi_{j_{1}})^{1-\frac{1}{p}}
f(\xi_{j_{2}})^{1-\frac{1}{p}}} ~h_{n}^{\frac{1}{p}
\{(d-\varepsilon)(\rho-p)-d \}} \big(1-2
h_{n}^{d-\varepsilon}\big)^{1-\frac{\rho}{p}}.
\end{equation*}
Finally, setting $k_{4}:=2 \left\lceil \frac{u_{0}}{\delta}
\right\rceil$ and $k_{5}:=\frac{4 p^{2} (2
a_{0})^{1/p}}{(\rho-p) \delta^{\rho/p}}$, one has
\begin{multline} \label{CovBound2}
n h_{n}^{d} ~\mathrm{Cov}\big(\hat{f}_{j_{1}},\hat{f}_{j_{2}}\big) \leq
-h_{n}^{d} f(\xi_{j_{1}}) f(\xi_{j_{2}}) + k_{4} \sqrt{f(\xi_{j_{1}}) f(\xi_{j_{2}}) (1-h_{n}^{d} f(\xi_{j_{1}})) (1-h_{n}^{d} f(\xi_{j_{2}}))} \\
+ 2 k\big(\dot{\xi}_{j_{1}}\big) \left( h_{n}^{\varepsilon} -
\left\lceil \frac{u_{0}}{\delta} \right\rceil h_{n}^{d} \right) +
k_{5} \sqrt{f(\xi_{j_{1}})^{1 - \frac{1}{p}} f(\xi_{j_{2}})^{1 -
\frac{1}{p}}} ~h_{n}^{\frac{1}{p} \{(d-\varepsilon)(\rho-p)-d
\}} \big(1-2 h_{n}^{d-\varepsilon}\big)^{1-\frac{\rho}{p}},
\end{multline}
which implies the desired result. The covariance is thus a
$O\big(1/\big(n h_{n}^{d}\big)\big)$ for any $\varepsilon$ in $\left[0,d \Big(1-\frac{1}{\rho-p}\Big)\right]$.
~~$\blacksquare$

\subsubsection{Proof of Theorem \ref{H-MISE1}} \noindent{Renewal sampling} -- By integrating over
$\pi_{nj}$ the right-hand side of (\ref{VarBound1}) and by
summing up over all hypercubes, we first derive an asymptotic
upper bound for IV. For some $\varepsilon \in \left[0, d
\Big(1-\frac{1}{\rho-p}\Big)\right]$,
\begin{multline*} \qquad n h_{n}^{d} \int_{\pi_{nj}} \mathrm{Var}\big(\hat{f}_{j} \big) \mathrm{d}x \leq h_{n}^{d}
\bigg\{ f(\xi_{j}) (1-h_{n}^{d} f(\xi_{j})) \{1 + k_{1}\} \\
+ k_{2} k\big(\dot{\xi}_{j}\big) h_{n}^{\varepsilon} +
k_{3} f(\xi_{j})^{1-\frac{1}{p}} ~h_{n}^{\frac{1}{p}
\{(d-\varepsilon)(\rho-p)-d\}} \bigg\}. \qquad
\end{multline*}
Then using the approximation of integral by Riemann sums, i.e.,
\begin{equation*} \sum_{j} h_{n}^{d}
f^{\kappa}(\xi_{j}) = \|f^{\kappa}\|_{1} + o(1), ~\kappa=1-\frac{1}{p}, 1, 2,  ~~\textrm{and}~~
\sum_{j} h_{n}^{d} k\big(\dot{\xi}_{j}\big) = \|k\|_{1}
+ o(1), \end{equation*} one has
\begin{equation} \label{IVbound1} n h_{n}^{d} ~\mathrm{IV}\big(\hat{f}_{n}^{H}\big) \leq \bigg\{1 + k_{1} + k_{2} \|k\|_{1} ~h_{n}^{\varepsilon}
+ k_{3} \big\|f^{1-\frac{1}{p}}\big\|_{1} ~h_{n}^{\frac{1}{p}
\{(d-\varepsilon)(\rho-p)-d\}} \bigg\} (1+o(1)).
\end{equation}
The two parts of the theorem follow from the choice $h_{n}=cn^{-1/(d+2)}$, $0<c<\infty$. So Lemma \ref{H-BCI} yields
\begin{equation*}
\lim_{n \rightarrow \infty} n^{\frac{2}{d+2}} ~\mathrm{ISB}\big(\hat{f}_{n}^{H}\big) = \frac{c^{2}}{12}
R_{d}\big(f'\big),
\end{equation*}
and combining with (\ref{IVbound1}), if $p=\rho-1 ~(\varepsilon=0)$,
we have
\begin{equation*}
\limsup_{n \rightarrow \infty} n^{\frac{2}{d+2}} ~\mathrm{MISE}\big(\hat{f}_{n}^{H}\big) \leq \frac{c^{2}}{12}
R_{d}\big(f'\big) + \frac{1}{c^{d}} \left\{1 + k_{1} + k_{2}
\|k\|_{1} + k_{3} \big\|f^{1-\frac{1}{p}}\big\|_{1}\right\}.
\end{equation*}
If $p < \rho-1 ~(\varepsilon > 0)$, we improve the
asymptotic constant:
\begin{equation*}
\limsup_{n \rightarrow \infty}
n^{\frac{2}{d+2}} ~\mathrm{MISE}\big(\hat{f}_{n}^{H}\big)
\leq \frac{c^{2}}{12} R_{d}\big(f'\big) + \frac{1}{c^{d}} \{1 +
k_{1}\}.
\end{equation*}
\newline
\noindent{Jittered sampling} -- Now, let us integrate over $\pi_{nj}$ the
right-hand side of (\ref{VarBound2}):
\begin{align*} n h_{n}^{d} \int_{\pi_{nj}} \mathrm{Var}\big(\hat{f}_{j}\big) \mathrm{d}x
\leq h_{n}^{d} \bigg\{& f(\xi_{j}) (1-h_{n}^{d} f(\xi_{j})) \{1
+ k_{4}\} + 2 k\big(\dot{\xi}_{j}\big) \left(
h_{n}^{\varepsilon} - \left\lceil \frac{u_{0}}{\delta} \right\rceil h_{n}^{d} \right) \nonumber \\
& +~ k_{5} f(\xi_{j})^{1-\frac{1}{p}} ~h_{n}^{\frac{1}{p}
\{(d-\varepsilon)(\rho-p)-d \}} \big(1-2
h_{n}^{d-\varepsilon}\big)^{1-\frac{\rho}{p}}\bigg\},
\end{align*}
for any $\varepsilon \in \left[0, d
\Big(1-\frac{1}{\rho-p}\Big)\right]$. Then sum up over all indexes
$j$ to obtain
\begin{equation*}
n h_{n}^{d} ~\mathrm{IV}\big(\hat{f}_{n}^{H}\big) \leq
\bigg\{1 + k_{4} + 2 \|k\|_{1} \left( h_{n}^{\varepsilon} -
\left\lceil \frac{u_{0}}{\delta} \right\rceil h_{n}^{d} \right) +
k_{5} \big\|f^{1-\frac{1}{p}}\big\|_{1} ~h_{n}^{\frac{1}{p}
\{(d-\varepsilon)(\rho-p)-d \}}\bigg\} (1+o(1)).
\end{equation*}
Therefore, if $p=\rho-1$, the
bin width choice $h_{n}=cn^{-1/(d+2)}$, $0<c<\infty$, entails
\begin{equation*}
\limsup_{n \rightarrow \infty} n^{\frac{2}{d+2}} ~\mathrm{MISE}\big(\hat{f}_{n}^{H}\big) \leq \frac{c^{2}}{12}
R_{d}\big(f'\big) + \frac{1}{c^{d}} \left\{1 + k_{4} + 2 \|k\|_{1} + k_{5} \big\|f^{1-\frac{1}{p}}\big\|_{1}\right\}.
\end{equation*}
If $p < \rho-1$, we get a better asymptotic constant:
\begin{equation*}
\limsup_{n \rightarrow \infty} n^{\frac{2}{d+2}} ~\mathrm{MISE}\big(\hat{f}_{n}^{H}\big) \leq \frac{c^{2}}{12}
R_{d}\big(f'\big) + \frac{1}{c^{d}} \{1 + k_{4}\}. ~~\blacksquare
\end{equation*}

\subsubsection{Variance bounds with high frequency sampling}
The period depends now on the sample size in that
$\delta_{n} \downarrow 0^{+}$ as $n \rightarrow \infty$. We start
by giving a new bound expression for the variance of
$\hat{f}_{n}^{H}(x)$ which depends upon $\gamma_{0}$.

\begin{lemma}[\textnormal{\bf high frequency sampling}] \label{H-MajVar3}
If $A_{0}(ii)$ and $A'_{1}(i)(ii)$ hold, then we obtain
\begin{multline} \label{VarBound3}
n h_{n}^{d} ~\mathrm{Var}\big(\hat{f}_{j}\big) \leq
f(\xi_{j}) (1 - h_{n}^{d} f(\xi_{j})) + 2
\varphi\big(\dot{\xi}_{j}\big) \Bigg(\sum_{r=1}^{r_{n}^{0}}
\frac{1}{r^{\gamma_{0}}}\Bigg) h_{n}^{d} \delta_{n}^{-\gamma_{0}}
+ \bigg\{2 u_{0} \|f\|_{\infty} f(\xi_{j}) \\ + 2
(u_{1}-u_{0}+\delta_{n}) k\big(\ddot{\xi}_{j}\big) \sup_{u \in
[u_{0},u_{1}]} \pi(u) + 2 k\big(\ddot{\xi}_{j}\big)
\int_{u_{1}}^{\infty} \pi(u) \mathrm{d}u
\bigg(1+\frac{\pi(u_{1})}{\int_{u_{1}}^{\infty} \pi(u) \mathrm{d}u}
\delta_{n}\bigg)\bigg\} h_{n}^{d} \delta_{n}^{-1},
\end{multline}
with $\big(\xi_{j},\dot{\xi}_{j},\ddot{\xi}_{j}\big) \in
\pi_{nj}^{3}$ and it entails that the variance is a
$O\big(1/\big(n h_{n}^{d}\big)\big)$ with the following choices $\delta_{n}^{*}(\gamma_{0})$ of
$\delta_{n}$:
\begin{equation*}
\delta_{n}^{*}(\gamma_{0})=d_{1} h_{n}^{d}
\mathds{1}_{\{\gamma_{0}<1\}} + d_{2} h_{n}^{d}
\ln\big(h_{n}^{-d}\big) \mathds{1}_{\{\gamma_{0}=1\}} + d_{3}
h_{n}^{d/\gamma_{0}} \mathds{1}_{\{\gamma_{0}>1\}}, ~0<d_{1},
d_{2}, d_{3}<\infty. \end{equation*}
\end{lemma}

\noindent \textbf{Proof of Lemma \ref{H-MajVar3}} ~The calculus of
$V_{n}$ remains identical. Now to upper bound $C_{n}$, we have to
make use of the local assumption $A'_{1}(i)$. Set
$r_{n}^{0}:=\lfloor u_{0} / \delta_{n} \rfloor$ and
$r_{n}^{1}:=\lfloor u_{1} / \delta_{n} \rfloor$, since $X^{T}$ is
stationary one may write
\begin{equation*} C_{n} = \frac{2}{h_{n}^{d}} \sum_{r=1}^{n-1} \left(1 - \frac{r}{n}\right)
\mathrm{Cov}\big(\mathds{1}_{\pi_{nj_{1}}}(X_{0}),\mathds{1}_{\pi_{nj_{2}}}(X_{r
\delta_{n}})\big) =: C_{n,1} + C_{n,2},
\end{equation*}
where \begin{align*} C_{n,1} & := \frac{2}{h_{n}^{d}}
\sum_{r=1}^{r_{n}^{0}} \left(1 - \frac{r}{n}\right)
\mathrm{Cov}\big(\mathds{1}_{\pi_{nj_{1}}}(X_{0}),\mathds{1}_{\pi_{nj_{2}}}(X_{r
\delta_{n}})\big), \\
C_{n,2} & := \frac{2}{h_{n}^{d}} \sum_{r=r_{n}^{0}+1}^{n-1}
\left(1 - \frac{r}{n}\right)
\mathrm{Cov}\big(\mathds{1}_{\pi_{nj_{1}}}(X_{0}),\mathds{1}_{\pi_{nj_{2}}}(X_{r
\delta_{n}})\big). \end{align*} First using $A'_{1}(i)$ we get
\begin{align*}
|C_{n,1}| & \leq \frac{2}{h_{n}^{d}} \sum_{r=1}^{r_{n}^{0}}
\iint_{\pi_{nj_{1}} \times \pi_{nj_{2}}} \bigg\{\sup_{y \in
\mathbb{R}^{d}} f_{r \delta_{n}}(x,y) + \|f\|_{\infty} f(x)\bigg\} \mathrm{d}x \mathrm{d}y \\
& \leq 2 \sum_{r=1}^{r_{n}^{0}} \int_{\pi_{nj_{1}}} \big\{\varphi(x) (r \delta_{n})^{-\gamma_{0}} + \|f\|_{\infty} f(x)\big\} \mathrm{d}x \\
& \leq 2 \varphi\big(\dot{\xi}_{j_{1}}\big)
\Bigg(\sum_{r=1}^{r_{n}^{0}} \frac{1}{r^{\gamma_{0}}}\Bigg)
h_{n}^{d} \delta_{n}^{-\gamma_{0}} + 2 u_{0} \|f\|_{\infty}
f(\xi_{j_{1}}) h_{n}^{d} \delta_{n}^{-1},
\end{align*}
where $\big(\xi_{j_{1}},\dot{\xi}_{j_{1}}\big) \in
\pi_{nj_{1}}^{2}$. Setting $k_{6}:=2 u_{0} \|f\|_{\infty}$, we
obtain
\begin{equation*} |C_{n,1}| \leq 2 \varphi\big(\dot{\xi}_{j_{1}}\big) \Bigg(
\sum_{r=1}^{r_{n}^{0}} \frac{1}{r^{\gamma_{0}}} \Bigg) h_{n}^{d}
\delta_{n}^{-\gamma_{0}} + k_{6} f(\xi_{j_{1}}) h_{n}^{d}
\delta_{n}^{-1}.
\end{equation*}
Then using $A'_{1}(ii)$
\begin{equation*}
|C_{n,2}| \leq \frac{2}{h_{n}^{d}}
\sum_{r=r_{n}^{0}+1}^{n-1} \iint_{\pi_{nj_{1}} \times \pi_{nj_{2}}} \sup_{y \in \mathbb{R}^{d}} |g_{r \delta_{n}}(x,y)| \mathrm{d}x \mathrm{d}y \leq 2 h_{n}^{d} k\big(\ddot{\xi}_{j_{1}}\big)
\Bigg[\sum_{r=r_{n}^{0}+1}^{r_{n}^{1}} \pi(r \delta_{n}) +
\sum_{r=r_{n}^{1}+1}^{n-1} \pi(r \delta_{n})\Bigg],
\end{equation*}
where $\ddot{\xi}_{j_{1}} \in \pi_{nj_{1}}$. On the one hand, one has
\begin{equation*} \sum_{r=r_{n}^{0}+1}^{r_{n}^{1}} \pi(r \delta_{n}) \leq \left(r_{n}^{1}-r_{n}^{0}\right) \sup_{u \in [u_{0},u_{1}]}
\pi(u) \leq (u_{1}-u_{0}) \sup_{u \in [u_{0},u_{1}]} \pi(u)
\bigg(1+\frac{\delta_{n}}{u_{1}-u_{0}}\bigg) \delta_{n}^{-1}.
\end{equation*}
On the other hand, the monotonicity of $\pi(\cdot)$ implies
\begin{equation*}
\sum_{r=r_{n}^{1}+1}^{n-1} \pi(r \delta_{n}) \leq \delta_{n}^{-1}
\sum_{r=r_{n}^{1}+1}^{n-1} \delta_{n} \pi(r \delta_{n}) \leq
\left\{ (u_{1} - r_{n}^{1} \delta_{n})
\pi\big(\left(r_{n}^{1}+1\right) \delta_{n}\big) +
\int_{u_{1}}^{\infty} \pi(u) \mathrm{d}u \right\} \delta_{n}^{-1}.
\end{equation*}
Setting $k_{7}:=2 (u_{1}-u_{0}) \sup_{u \in [u_{0},u_{1}]} \pi(u)$
and $k_{8}:=2 \int_{u_{1}}^{\infty} \pi(u) \mathrm{d}u$, we thus obtain
\begin{equation*} |C_{n,2}| \leq k_{7}
k\big(\ddot{\xi}_{j_{1}}\big)
\bigg(1+\frac{\delta_{n}}{u_{1}-u_{0}}\bigg) h_{n}^{d}
\delta_{n}^{-1} + k_{8} k\big(\ddot{\xi}_{j_{1}}\big) \bigg(1 +
\frac{\pi(u_{1})}{\int_{u_{1}}^{\infty} \pi(u) \mathrm{d}u} \delta_{n}
\bigg) h_{n}^{d} \delta_{n}^{-1}.
\end{equation*}
Thence
\begin{multline} \label{CovBound3}
n h_{n}^{d} ~\mathrm{Cov}\big(\hat{f}_{j_{1}},\hat{f}_{j_{2}}\big) \leq
-h_{n}^{d} f(\xi_{j_{1}}) f(\xi_{j_{2}}) + 2
\varphi\big(\dot{\xi}_{j_{1}}\big) \Bigg(\sum_{r=1}^{r_{n}^{0}}
\frac{1}{r^{\gamma_{0}}}\Bigg) h_{n}^{d} \delta_{n}^{-\gamma_{0}}
+ k_{6} f(\xi_{j_{1}}) h_{n}^{d} \delta_{n}^{-1} \\
+ k_{7} k\big(\ddot{\xi}_{j_{1}}\big)
\bigg(1+\frac{\delta_{n}}{u_{1}-u_{0}}\bigg) h_{n}^{d}
\delta_{n}^{-1} + k_{8} k\big(\ddot{\xi}_{j_{1}}\big) \bigg(1 +
\frac{\pi(u_{1})}{\int_{u_{1}}^{\infty} \pi(u) \mathrm{d}u}
\delta_{n}\bigg) h_{n}^{d} \delta_{n}^{-1},
\end{multline}
which leads to the desired result. Using (\ref{CovBound3}), we
also deduce the optimal choices $\delta_{n}^{*}(\gamma_{0})$ of
$\delta_{n}$ i.e. the smallest values of $\delta_{n}$ so
that $C_{n}$ is a $O(1)$. These choices are given by
(\ref{DeltaOpt}) in accordance with the values of $\gamma_{0}$.
~~$\blacksquare$

\subsubsection{Proof of Theorem \ref{H-MISE2Opt}} By integrating over
$\pi_{nj}$ the right-hand side of (\ref{VarBound3}):
\begin{multline*} n h_{n}^{d} \int_{\pi_{nj}} \mathrm{Var}\big(\hat{f}_{j}\big) \mathrm{d}x
\leq h_{n}^{d} \Bigg\{f(\xi_{j}) (1-h_{n}^{d} f(\xi_{j})) + 2
\varphi(\dot{\xi}_{j}) \Bigg(\sum_{r=1}^{r_{n}^{0}}
\frac{1}{r^{\gamma_{0}}}\Bigg) h_{n}^{d} \delta_{n}^{-\gamma_{0}} + k_{6} f(\xi_{j}) h_{n}^{d} \delta_{n}^{-1} \\
+ k_{7} k(\ddot{\xi}_{j})
\bigg(1+\frac{\delta_{n}}{u_{1}-u_{0}}\bigg) h_{n}^{d}
\delta_{n}^{-1} + k_{8} k(\ddot{\xi}_{j})
\bigg(1+\frac{\pi(u_{1})}{\int_{u_{1}}^{\infty} \pi(u) \mathrm{d}u}
\delta_{n}\bigg) h_{n}^{d} \delta_{n}^{-1}\Bigg\}.
\end{multline*}
Then let us sum up over all indexes $j$. Since $\varphi$ is
Riemann-integrable, we obtain
\begin{equation*}
n h_{n}^{d} ~\mathrm{IV}\big(\hat{f}_{n}^{H}\big) \leq
\Bigg\{1 + 2 \|\varphi\|_{1} \Bigg(\sum_{r=1}^{r_{n}^{0}}
\frac{1}{r^{\gamma_{0}}}\Bigg) h_{n}^{d} \delta_{n}^{-\gamma_{0}}
+ k_{6} h_{n}^{d} \delta_{n}^{-1} + (k_{7}+k_{8}) \|k\|_{1}
h_{n}^{d} \delta_{n}^{-1}\Bigg\} (1+o(1)).
\end{equation*}
According to the values of $\gamma_{0}$, we derive all asymptotic
bounds with optimal choices of $\delta_{n}$:
\begin{itemize}
\item[--] if $\gamma_{0}<1$, the choice
$\delta_{n}^{*}(\gamma_{0})=d_{1} h_{n}^{d}$, $0 < d_{1} <
\infty$, entails
$$\Bigg(\sum_{r=1}^{r_{n}^{0}} \frac{1}{r^{\gamma_{0}}}\Bigg)
h_{n}^{d} \delta_{n}^{-\gamma_{0}} \leq \frac{1}{d_{1}}
\frac{u_{0}^{1-\gamma_{0}}}{1-\gamma_{0}} \bigg[1-
\frac{d_{1}^{1-\gamma_{0}} \gamma_{0}}{u_{0}^{1-\gamma_{0}}}
h_{n}^{d(1-\gamma_{0})}\bigg];$$ \item[--] if $\gamma_{0}=1$, the
choice $\delta_{n}^{*}(\gamma_{0})=d_{2} h_{n}^{d}
\ln\left(h_{n}^{-d}\right)$, $0 < d_{2} < \infty$, entails
$$\Bigg(\sum_{r=1}^{r_{n}^{0}} \frac{1}{r}\Bigg) h_{n}^{d}
\delta_{n}^{-1} \leq \frac{1}{d_{2}} \bigg[1 -
\ln\big(h_{n}^{-d}\big) \ln\left(e\frac{u_{0}}{d_{2}}
\ln\big(h_{n}^{-d}\big)\right)\bigg];$$ \item[--] if
$\gamma_{0}>1$, the choice $\delta_{n}^{*}(\gamma_{0})=d_{3}
h_{n}^{d/\gamma_{0}}$, $0 < d_{3} < \infty$, entails
$$\Bigg(\sum_{r=1}^{r_{n}^{0}} \frac{1}{r^{\gamma_{0}}}\Bigg)
h_{n}^{d} \delta_{n}^{-\gamma_{0}} \leq
\frac{\gamma_{0}}{d_{3}^{\gamma_{0}} (\gamma_{0}-1)}
\bigg[1-\frac{1}{\gamma_{0} u_{0}^{\gamma_{0}-1}}
\delta_{n}^{\gamma_{0}-1}\bigg].$$
\end{itemize}
So setting
$$C_{\gamma_{0}}:=\frac{1}{d_{1}}\bigg\{ \frac{2 \|\varphi\|_{1} u_{0}^{1-\gamma_{0}}}{1-\gamma_{0}} + k_{6} + (k_{7}+k_{8}) \|k\|_{1}
\bigg\} \mathds{1}_{\{\gamma_{0}<1\}} + \frac{2
\|\varphi\|_{1}}{d_{2}} \mathds{1}_{\{\gamma_{0}=1\}} + \frac{2
\|\varphi\|_{1} \gamma_{0}}{d_{3}^{\gamma_{0}} (\gamma_{0}-1)}
\mathds{1}_{\{\gamma_{0}>1\}},$$ it remains to choose $h_{n}=cn^{-1/(d+2)}$, $0<c<\infty$, as
in Theorem \ref{H-MISE1}:
\begin{equation*}
\limsup_{n \rightarrow \infty} n^{\frac{2}{d+2}} ~\mathrm{MISE}\big(\hat{f}_{n}^{H}\big) \leq \frac{c^{2}}{12}
R_{d}\big(f'\big) + \frac{1}{c^{d}}\big(1+C_{\gamma_{0}}\big),
\end{equation*}
and we can also improve our asymptotic constant for any choice of
$\delta_{n}$ such that $\delta_{n}/\delta_{n}^{*}(\gamma_{0}) \rightarrow \infty$ as $n \rightarrow \infty$:
\begin{equation*}
\limsup_{n \rightarrow \infty} n^{\frac{2}{d+2}}
~\mathrm{MISE}\big(\hat{f}_{n}^{H}\big) \leq \frac{c^{2}}{12}
R_{d}\big(f'\big) + \frac{1}{c^{d}}. ~~\blacksquare
\end{equation*}

\subsection{Frequency polygon}

\subsubsection{Proof of Theorem \ref{FP-MISE1}}
\noindent Renewal sampling -- First observe that
\begin{equation*}
\int_{\mathbb{R}} \mathrm{Var}\big(\hat{f}_{n}^{FP}(x)\big) \mathrm{d}x =
\sum_{j} \int_{[c_{j},c_{j+1}[}
\mathrm{Var}\big(\hat{f}_{n}^{FP}(x)\big) \mathrm{d}x,
\end{equation*}
where
\begin{multline*}
\int_{c_{j}}^{c_{j+1}} \mathrm{Var}\big(\hat{f}_{n}^{FP}(x)\big)
\mathrm{d}x = \frac{1}{h_{n}^{2}} \int_{c_{j}}^{c_{j+1}} \bigg\{
(x-c_{j})^{2} \mathrm{Var}\big(\hat{f}_{j+1}\big) + (c_{j+1}-x)^{2} \mathrm{Var}\big(\hat{f}_{j}\big) \\
+ 2 (x-c_{j})(c_{j+1}-x)
\mathrm{Cov}\big(\hat{f}_{j},\hat{f}_{j+1}\big) \bigg\} \mathrm{d}x.
\end{multline*}
For any $j \in \mathbb{Z}$, let us denote by $\overline{V}_{j}$
(respectively $\overline{C}_{j,j+1}$) an upper bound expression
for $n h_{n} \mathrm{Var}\big(\hat{f}_{j}\big)$
(respectively $n h_{n} \mathrm{Cov}\big(\hat{f}_{j},\hat{f}_{j+1}\big)$) that is
independent of $x$. We get
\begin{equation} \label{FPVarBound}
n h_{n} \int_{c_{j}}^{c_{j+1}}
\mathrm{Var}\big(\hat{f}_{n}^{FP}(x)\big) \mathrm{d}x \leq \frac{h_{n}}{3}
\left\{\overline{V}_{j} + \overline{V}_{j+1} +
\overline{C}_{j,j+1}\right\}.
\end{equation}
Insert now both expressions (\ref{VarBound1}) and
(\ref{CovBound1}) in (\ref{FPVarBound}), then for $\varepsilon
\in \left[0,1-\frac{1}{\rho-p}\right]$,
\begin{align*}
& n h_{n} \int_{c_{j}}^{c_{j+1}} \mathrm{Var}\big(\hat{f}_{n}^{FP}(x)\big) \mathrm{d}x \\
& \leq \frac{h_{n}}{3} \bigg\{ f(\xi_{j}) (1-h_{n} f(\xi_{j}))
\{1 + k_{1}\} + k_{2} k(\dot{\xi}_{j}) h_{n}^{\varepsilon} + k_{3}
f(\xi_{j})^{1-\frac{1}{p}} ~h_{n}^{\frac{1}{p} \{
(1-\varepsilon)(\rho-p)-1\}} \bigg\} \\
& ~+ \frac{h_{n}}{3} \bigg\{ f(\xi_{j+1}) (1-h_{n} f(\xi_{j+1}))
\{1 + k_{1}\} + k_{2} k(\dot{\xi}_{j+1}) h_{n}^{\varepsilon} +
k_{3} f(\xi_{j+1})^{1-\frac{1}{p}} ~h_{n}^{\frac{1}{p} \{
(1-\varepsilon)(\rho-p)-1 \}}
\bigg\} \\
& ~+ \frac{h_{n}}{3} \bigg\{-h_{n} f(\xi_{j}) f(\xi_{j+1}) +
k_{1} \sqrt{f(\xi_{j}) f(\xi_{j+1}) (1-h_{n} f(\xi_{j})) (1-h_{n} f(\xi_{j+1}))} \\
& \quad \quad \quad + k_{2} k(\dot{\xi}_{j}) h_{n}^{\varepsilon} +
k_{3} \sqrt{f(\xi_{j})^{1-\frac{1}{p}} f(\xi_{j+1})^{1 -
\frac{1}{p}}} ~h_{n}^{\frac{1}{p} \{(1-\varepsilon)(\rho-p)-1
\}}\bigg\}.
\end{align*}
We bound the IV of $\hat{f}_{n}^{FP}$ by summing up over all indexes $j$. So for $\varepsilon \in
\left[0,1-\frac{1}{\rho-p}\right]$,
\begin{equation} \label{FP-IVbound1}
n h_{n} ~\mathrm{IV}\big(\hat{f}_{n}^{FP}\big) \leq \bigg\{
\frac{2}{3} + k_{1} + k_{2} \|k\|_{1} ~h_{n}^{\varepsilon} + k_{3}
\big\|f^{1-\frac{1}{p}}\big\|_{1} ~h_{n}^{\frac{1}{p}
\{(1-\varepsilon)(\rho-p)-1\}} \bigg\} (1+o(1)).
\end{equation}
Now the bin width choice $h_{n}=cn^{-1/5}$, $0<c<\infty$, in Lemma \ref{FP-BCI} yields first
\begin{equation*}
\lim_{n \rightarrow \infty} n^{\frac{4}{5}} ~\mathrm{ISB}\big(\hat{f}_{n}^{FP}\big) = \frac{49}{2880} c^{4}
R\big(f''\big),
\end{equation*}
and combining with (\ref{FP-IVbound1}), if $p = \rho-1 ~(\varepsilon=0)$:
\begin{equation*}
\limsup_{n \rightarrow \infty} n^{\frac{4}{5}} ~\mathrm{MISE}\big(\hat{f}_{n}^{FP}\big) \leq \frac{49}{2880} c^{4}
R\big(f''\big) + \frac{1}{c} \left\{ \frac{2}{3} + k_{1} + k_{2} \|k\|_{1} +
k_{3} \big\|f^{1-\frac{1}{p}}\big\|_{1} \right\}.
\end{equation*}
If $p < \rho-1 ~(0 < \varepsilon < 1)$:
\begin{equation*}
\limsup_{n \rightarrow \infty} n^{\frac{4}{5}} ~\mathrm{MISE}\big(\hat{f}_{n}^{FP}\big) \leq \frac{49}{2880} c^{4}
R\big(f''\big) + \frac{1}{c} \left\{ \frac{2}{3} + k_{1} \right\}.
\end{equation*}
\newline
\noindent Jittered sampling -- The outlines of the proof are unchanged.
Insert both expressions (\ref{VarBound2}) and (\ref{CovBound2})
in (\ref{FPVarBound}) and sum up over all indexes $j$, then
it follows that for $\varepsilon \in \left[0,1-\frac{1}{\rho-p}\right]$,
\begin{multline*} n h_{n} ~\mathrm{IV}\big(\hat{f}_{n}^{FP}\big) \leq
\bigg\{ \frac{2}{3} + k_{4} + 2 \|k\|_{1} \left( h_{n}^{\varepsilon} -
\left\lceil \frac{u_{0}}{\delta} \right\rceil  h_{n} \right) \\ +
k_{5} \big\|f^{1-\frac{1}{p}}\big\|_{1} ~h_{n}^{\frac{1}{p}
\{(1-\varepsilon)(\rho-p)-1 \}} \left(1-2 h_{n}^{1-\varepsilon}\right)^{1-\frac{\rho}{p}} \bigg\} (1+o(1)). \end{multline*}
Take $h_{n}=cn^{-1/5}$, $0 < c < \infty$, then if $p=\rho-1$:
\begin{equation*}
\limsup_{n \rightarrow \infty} n^{\frac{4}{5}} ~\mathrm{MISE}\big(\hat{f}_{n}^{FP}\big) \leq \frac{49}{2880} c^{4}
R\big(f''\big) + \frac{1}{c} \left\{ \frac{2}{3} + k_{4} + 2 \|k\|_{1} +
k_{5} \big\|f^{1-\frac{1}{p}}\big\|_{1} \right\}.
\end{equation*}
If $p < \rho-1$:
\begin{equation*}
\limsup_{n \rightarrow \infty} n^{\frac{4}{5}} ~\mathrm{MISE}\big(\hat{f}_{n}^{FP}\big) \leq \frac{49}{2880} c^{4}
R\big(f''\big) + \frac{1}{c} \left\{ \frac{2}{3} + k_{4} \right\}.
~~\blacksquare
\end{equation*}

\subsubsection{Proof of Theorem \ref{FP-MISE2Opt}}
Insert now both expressions (\ref{VarBound3}) and
(\ref{CovBound3}) in (\ref{FPVarBound}) and sum up over all indexes $j$, we get
\begin{equation*} n h_{n} ~\mathrm{IV}\big(\hat{f}_{n}^{FP}\big) \leq \Bigg\{\frac{2}{3} + 2 \|\varphi\|_{1}
\Bigg(\sum_{r=1}^{r_{n}^{0}} \frac{1}{r^{\gamma_{0}}}\Bigg) h_{n}
\delta_{n}^{-\gamma_{0}} + k_{6} h_{n} \delta_{n}^{-1} +
(k_{7}+k_{8}) \|k\|_{1} h_{n} \delta_{n}^{-1}\Bigg\} (1+o(1)).
\end{equation*}
Then $h_{n}=cn^{-1/5}$, $0<c<\infty$, together with
the optimal choices $\delta_{n}^{*}$ of $\delta_{n}$ yield
\begin{equation*}
\limsup_{n \rightarrow \infty} n^{\frac{4}{5}} ~\mathrm{MISE}\big(\hat{f}_{n}^{FP}\big) \leq \frac{49}{2880} c^{4}
R\big(f''\big) + \frac{1}{c} \bigg\{\frac{2}{3} + C_{\gamma_{0}}\bigg\},
\end{equation*}
and, if $\delta_{n}$ is such that $\delta_{n}/\delta_{n}^{*}(\gamma_{0}) \rightarrow \infty$ as $n \rightarrow \infty$, $C_{\gamma_{0}}$ is removable in the limiting bound. ~~$\blacksquare$

\bibliographystyle{plain}
\bibliography{RefProjSamp}

\end{document}